\documentclass[12pt]{article}
\usepackage{color}
\definecolor{darkblue}{rgb}{0.00,0.25,0.50}
\usepackage[colorlinks,filecolor=blue,citecolor=darkblue]{hyperref}

\usepackage{amsfonts,amssymb,amsmath,amsthm}
\usepackage{url}
\usepackage{enumerate}
\usepackage[ukrainian, russian, english]{babel}
\usepackage[cp1251]{inputenc}
\begin{document} \selectlanguage{ukrainian}
\thispagestyle{empty}

\title{}

UDC 517.51 \vskip 5mm

\begin{center}
\textbf{\Large Order estimates of approximation characteristics of functions from anisotropic Nikol'skii--Besov classes }
\end{center}

\vskip 3mm

\begin{center}
\textbf{\Large  Порядкові оцінки апроксимативних характеристик функцій з
анізотропних класів  Нікольського--Бєсова}
\end{center}
\vskip0.5cm

\begin{center}
 S.~Ya.~Yanchenko\\ \emph{\small
Institute of Mathematics of NAS of
Ukraine, Kyiv}
\end{center}
\begin{center}
C.~Я.~Янченко \\
\emph{\small Інститут математики НАН України, Київ}
\end{center}
\vskip0.5cm

\begin{abstract}

We obtained exact order estimates of the deviation of functions from anisotropic Nikol'skii--Besov classes $B^{\boldsymbol{r}}_{p,\theta}(\mathbb{R}^d)$ from their
sections of the Fourier integral. The error of the approximation is estimated in the metric of
Lebesgue space  $L_{\infty}(\mathbb{R}^d)$.
\vskip 3 mm

Одержано точні за порядком оцінки відхилення функцій з анізотропних класів
Нікольського--Бєсова $B^{\boldsymbol{r}}_{p,\theta}(\mathbb{R}^d)$ від їх відрізків інтеграла Фур'є.  Похибка наближення вимірюється у метриці простору  $L_{\infty}(\mathbb{R}^d)$.
\end{abstract}

\vskip 0.5 cm

%%%%%%%%%%%%%%%%%%%%%%%%%%%%%%%%%%%%%%%%%%%%%%%%%%%%%%%%%%%%%%%%%%%%%%%%%

\textbf{Вступ.}  У роботі встановлюються  порядкові оцінки деяких апроксимативних характеристик класів функцій багатьох змінних з анізотропних просторів Нікольського--Бєсова $B^{\boldsymbol{r}}_{p,\theta}(\mathbb{R}^d)$,  де $\boldsymbol{r}=(r_1,\ldots,r_d)$, $r_j>0$, ${j=\overline{1,d}}$, $1\leqslant p, \ \theta\leqslant\infty$. Простори $B^{\boldsymbol{r}}_{p,\theta}(\mathbb{R}^d)$ були введені О.\,В.~Бєсовим~\cite{Besov_1961} і ${B^{\boldsymbol{r}}_{p,\infty}(\mathbb{R}^d)=H^{\boldsymbol{r}}_p(\mathbb{R}^d)}$, де
$H^{\boldsymbol{r}}_p(\mathbb{R}^d)$~--- простори, які ввів
 С.\,М.~Нікольський~\cite{Nikolsky_1951}.

 Вихідні означення просторів $H^{\boldsymbol{r}}_p(\mathbb{R}^d)$ та $B^{\boldsymbol{r}}_{p,\theta}(\mathbb{R}^d)$
 у згаданих роботах були дані в термінах певних обмежень на  модулі гладкості функцій з цих просторів. У подальших дослідженнях нам  буде зручно користуватися еквівалентним означення просторів $B^{\boldsymbol{r}}_{p,\theta}(\mathbb{R}^d)$, яке було встановлене  П.\,І.~Лізоркіним~\cite{Lizorkin_1968_sib}, що базується на  використані перетворення Фур'є.

Зауважимо, що анізотропні простори Нікольського--Бєсова функцій багатьох змінних, що визначені на $\mathbb{R}^d$ з точки зору знаходження точних за порядком значень деяких апроксимативних характеристик досліджувалися, зокрема, у роботах~\cite{Jiang_Yanjie_Liu_Yongping_JAT_2000}, \cite{Jiang_Yanjie_AM_2002}, а у випадку $r_1=\ldots=r_d=r$, тобто ізотропні простори Нікольського--Бєсова $B^r_{p,\theta}(\mathbb{R}^d)$, досліджувалися  у роботах \cite{Yanchenko_Zb_2010TF}, \cite{Yanchenko_2015UMG}.

\textbf{1. Основні позначення та означення класів
Нікольського--Бєсова.} Нехай  $\mathbb{R}^d$, $d\geqslant 1$~---
$d$-вимірний евклідів простір з елементами
${\boldsymbol{x}=(x_1,...,x_d)}$,
${(\boldsymbol{x},\boldsymbol{y})=x_1y_1+...+x_dy_d}$.
${L_p(\mathbb{R}^d)}$, ${1\leqslant p\leqslant\infty}$,~--- простір
вимірних на $\mathbb{R}^d$ функцій ${f(\boldsymbol{x})=f(x_1,...,x_d)}$ зі скінченною
нормою
 $$
\|f\|_{L_p(\mathbb{R}^d)}=\|f\|_p:=
\left(\int\limits_{\mathbb{R}^{d}}|f(\boldsymbol{x})|^{p}d\boldsymbol{x}
\right) ^{\frac{1}{p}}, \ 1\leqslant p<\infty,
 $$
 $$
  \|f\|_{L_\infty(\mathbb{R}^d)}=\|f\|_{\infty}:=\mathop {\rm ess \sup}\limits_{\boldsymbol{x}\in \mathbb{R}^d}
  |f(\boldsymbol{x})|.
 $$

Нехай $S=S(\mathbb{R}^d)$~--- простір Л.~Шварца основних нескінченно
диференційовних на $\mathbb{R}^d$ комплекснозначних функцій
$\varphi$, що спадають на нескінченності разом зі своїми похідними
швидше за будь-який степінь функції $(x_1^2+\ldots+x_d^2)^{-\frac{1}{2}}$ (див., наприклад,
\cite{Lizorkin_69}). Через $S'$ позначимо
простір лінійних неперервних функціоналів на $S$. Зазначимо, що
елементами простору $S'$ є узагальнені функції. Якщо $f\in S'$ і
$\varphi \in S$, то $\langle f,\varphi\rangle$ позначає значення $f$
на $\varphi$.

Перетворення Фур'є $\mathfrak{F}\varphi\colon S\rightarrow S$
визначається згідно з формулою
$$
(\mathfrak{F}\varphi)(\boldsymbol{\lambda})=\frac{1}{(2\pi)^{d/2}}\int
 \limits_{\mathbb{R}^d}\varphi(\boldsymbol{t})
 e^{-i(\boldsymbol{\lambda},\boldsymbol{t})}d\boldsymbol{t}
 \equiv \widetilde{\varphi}(\boldsymbol{\lambda}).
$$

Обернене перетворення Фур'є $\mathfrak{F}^{-1}\varphi\colon S\rightarrow S$ задається таким чином:
$$
(\mathfrak{F}^{-1}\varphi)(\boldsymbol{t})=\frac{1}{(2\pi)^{d/2}}
\int \limits_{\mathbb{R}^d}\varphi(\boldsymbol{\lambda})
e^{i(\boldsymbol{\lambda},\boldsymbol{t})}d\boldsymbol{\lambda}\equiv
\widehat{\varphi}(\boldsymbol{t}).
$$

Перетворення Фур'є узагальнених функцій $f\in S'$ (для нього ми
зберігаємо те ж позначення)  визначається згідно з формулою
$$
\langle \mathfrak{F}f,\varphi\rangle=\langle f,\mathfrak{F}\varphi
\rangle  \ \ \ (\langle \widetilde{f},\varphi\rangle=\langle
f,\widetilde{\varphi} \rangle),
$$
де $\varphi \in S$.

Обернене перетворення Фур'є узагальненої функції  $f\in S'$ також позначимо
$\mathfrak{F}^{-1}f$,  і визначається воно аналогічно до прямого
перетворення Фур'є згідно з формулою
$$
\langle \mathfrak{F}^{-1}f,\varphi\rangle=\langle
f,\mathfrak{F}^{-1}\varphi \rangle  \ \ \ (\langle
\widehat{f},\varphi\rangle=\langle f,\widehat{\varphi} \rangle).
$$

Зазначимо, що кожна функція $f \in L_p(\mathbb{R}^d)$, $1\leqslant p \leqslant \infty$,
 визначає лінійний неперервний функціонал на $S$ згідно з формулою
$$
\langle f,\varphi \rangle = \int \limits_{\mathbb{R}^d}
 f(\boldsymbol{x})\varphi(\boldsymbol{x}) d\boldsymbol{x}, \ \ \varphi\in S,
$$
і, як наслідок, у цьому сенсі вона є елементом $S'$. Тому
перетворення Фур'є функції $f \in L_p(\mathbb{R}^d)$, $1\leqslant p \leqslant \infty$, можна
розглядати як перетворення Фур'є узагальненої функції $\langle f,\varphi \rangle$.

Носієм узагальненої функції $f$ будемо називати замикання
$\overline{\mathfrak{N}}$ такої множини точок
$\mathfrak{N}\subset\mathbb{R}^d$, що для довільної $\varphi \in S$,
яка дорівнює нулю в $\overline{\mathfrak{N}}$, виконується рівність
$\langle f,\varphi \rangle = 0$. Носій узагальненої функції $f$
будемо позначати через $\mbox{supp}\, f$. Також будемо говорити, що функція $f$
 зосереджена на множині $G$, якщо $\mbox{supp}\, f \subseteq G$.

У подальшому будемо користуватися такими позначеннями. Нехай функція $f$ представлена інтегралом Фур'є
$$
f(\boldsymbol{x})=\frac{1}{(2\pi)^{d/2}}
\int\limits_{\mathbb{R}^d}\tilde{f}(\boldsymbol{\lambda})e^{i(\boldsymbol{\lambda}, \boldsymbol{x})}d\boldsymbol{\lambda}.
$$
Тоді відрізком інтегралу Фур'є функції $f$  назвемо вираз
$$
S_{\boldsymbol{\sigma}}(f)=\frac{1}{(2\pi)^{d/2}}\int\limits_{-\sigma_1}^{\sigma_1}\ldots
\int\limits_{-\sigma_d}^{\sigma_d}
\tilde{f}(\boldsymbol{\lambda})e^{i(\boldsymbol{\lambda}, \boldsymbol{x})}d\boldsymbol{\lambda},
$$
де $\tilde{f}(\boldsymbol{\lambda})$~--- перетворення Фур'є функції $f\in L_p(\mathbb{R}^d)$.

Крім того, для $S_{\boldsymbol{\sigma}}(f)$ можемо записати (див. \cite{Lizorkin_69})
$$
S_{\boldsymbol{\sigma}}(f)=\frac{1}{\pi^{d}}\int\limits_{\mathbb{R}^d} f(\boldsymbol{y}) \prod\limits_{j=1}^d \frac{\sin{\sigma_j(x_j-y_j)}}{x_j-y_j} d\boldsymbol{y}.
$$
Таким чином $S_{\boldsymbol{\sigma}}(f)$~--- ціла функція степеня $\boldsymbol{\sigma}$.

Нехай $D_{\boldsymbol{a}^s}=D_{a_1^{s},\ldots, a_d^{s}}$~--- паралелепіпед:
$|\lambda_j|<a_j^s$, $j=\overline{1,d}$, $s\in \mathbb{N}\cup\{0\}$, а ${\Gamma_{\boldsymbol{a}^s}=D_{\boldsymbol{a}^s}-D_{\boldsymbol{a}^{s-1}}}$ при $s \geqslant 1$ і $\Gamma_{\boldsymbol{a}^0}=D_{\boldsymbol{a}^0}$. Покладемо
$$
f_{\boldsymbol{a}^s}=S_{\boldsymbol{a}^s}(f)-S_{\boldsymbol{a}^{s-1}}(f)= \int\limits_{\Gamma_{\boldsymbol{a}^s}}\tilde{f}(\boldsymbol{\lambda})e^{i(\boldsymbol{\lambda}, \boldsymbol{x})}d\boldsymbol{\lambda}, s\geqslant 1,
$$
і
$$
f_{\boldsymbol{a}^0}=S_{\boldsymbol{a}^0}(f)= \int\limits_{\Gamma_{\boldsymbol{a}^0}}\tilde{f}(\boldsymbol{\lambda})e^{i(\boldsymbol{\lambda}, \boldsymbol{x})}d\boldsymbol{\lambda}.
$$
Представлення функції $f$ у вигляді
$$
f=f_{\boldsymbol{a}^0}+\sum\limits_{s=1}^{\infty}f_{\boldsymbol{a}^s}= \sum\limits_{s=0}^{\infty}f_{\boldsymbol{a}^s}
$$
будемо називати розшаруванням $f$  ($\boldsymbol{a}$-розшаруванням $f$). У випадку, коли $f\in L_p$, ${p>2}$, $S_{\boldsymbol{a}^s}(f)$ розуміють, взагалі кажучи, як результат дії на $f$ деякого оператора, який в образах Фур'є зводиться до множення на характеристичну функцію області $D_{\boldsymbol{a}^s}$ (див. \cite{Lizorkin_69}, (\S3, гл.1)).

Далі для вектора $\boldsymbol{r}=(r_1,\ldots,r_d)$, $r_j>0$, $j=\overline{1,d}$, введемо величину
\begin{equation}\label{g(r)}
g(\boldsymbol{r})=\Bigg(\frac{1}{d}\sum\limits_{j=1}^{d}\frac{1}{r_j}\Bigg)^{-1}.
\end{equation}
Зауважимо, що при $r_1=r_2=\ldots=r_d=r$ маємо $g(\boldsymbol{r})=r$.

Тоді  анізотропні простори $B^{\boldsymbol{r}}_{p,\theta}(\mathbb{R}^d)$ можна означити таким чином~\cite{Lizorkin_1968_sib}:
$$
B^{\boldsymbol{r}}_{p,\theta}(\mathbb{R}^d)=\left\{ f\in L_p(\mathbb{R}^d)\colon \|f\|_{B^{\boldsymbol{r}}_{p,\theta}}<\infty \right\},
$$
де
\begin{equation}\label{f-norm-dek}
 \|f\|_{B^{\boldsymbol{r}}_{p,\theta}(\mathbb{R}^d)}\asymp
 \left(\sum\limits_{s=0}^{\infty}b^{s\theta}\|f_{\boldsymbol{a}^s}\|_p^{\theta}
 \right)^{\frac{1}{\theta}}<\infty,   \ \ \mbox{при} \ \ 1\leqslant \theta <\infty,
\end{equation}
\begin{equation}\label{f-norm-dek-inf}
 \|f\|_{B^{\boldsymbol{r}}_{p,\infty}(\mathbb{R}^d)}\asymp
 \sup\limits_{s\geqslant 0} b^{s}\|f_{\boldsymbol{a}^s}\|_p <\infty,
\end{equation}
а $b=2^{g(\boldsymbol{r})}$ і $a_j=2^{g(\boldsymbol{r})/r_j}$, $j=\overline{1,d}$.

Далі, зберігаючи ті самі позначення, будемо розглядати класи $B^{\boldsymbol{r}}_{p,\theta}(\mathbb{R}^d)$, тобто одиничні кулі у просторах $B^{\boldsymbol{r}}_{p,\theta}(\mathbb{R}^d)$:
$$
B^{\boldsymbol{r}}_{p,\theta}(\mathbb{R}^d):= \left\{f\in L_p(\mathbb{R}^d)\colon
\|f\|_{B^{\boldsymbol{r}}_{p,\theta}(\mathbb{R}^d)}\leqslant 1\right\}.
$$
Крім цього, для спрощення записів, замість $B^{\boldsymbol{r}}_{p,\theta}(\mathbb{R}^d)$ та
$H^{\boldsymbol{r}}_{p}(\mathbb{R}^d)$ будемо  використовувати  позначення
$B^{\boldsymbol{r}}_{p,\theta}$ та $H^{\boldsymbol{r}}_{p}$.

Зазначимо, що при встановленні результатів важливим є те, що простори $B^{\boldsymbol{r}}_{p,\theta}$ зі зростанням значення параметра $\theta$
розширюються (див., наприклад, \cite[с.~278]{Nikolsky_1969_book}),
тобто
\begin{equation}\label{vklad}
 B^{\boldsymbol{r}}_{p,1}\subset B^{\boldsymbol{r}}_{p,\theta}\subset B^{\boldsymbol{r}}_{p,\theta'}\subset
 B^{\boldsymbol{r}}_{p,\infty}=H^{\boldsymbol{r}}_p, \ \ \ 1\leqslant\theta<\theta' \leqslant
 \infty.
\end{equation}

\vskip 3 mm
\textbf{2. Допоміжні твердження та основний результат.}

Важливе значення при доведенні одержаного результату відіграє теорема встановлена О.\,В.~Бєсовим~\cite{Besov_1961} (теорема 2.1), яку сформулюємо у такій формі.

\vskip 1 mm
\bf Теорема А. \it  Нехай
$1\leqslant p \leqslant p' \leqslant \infty$, $\theta'\geqslant \theta$,
$$
\kappa=1 - \left(\frac{1}{p}-\frac{1}{p'}\right)\sum\limits_{j=1}^{d}\frac{1}{r_j}>0.
$$
Тоді, якщо $f\in B^{\boldsymbol{r}}_{p,\theta}$, то $f\in B^{\boldsymbol{\rho}}_{p',\theta'}$, де $\rho_j=r_j\kappa$, $j=\overline{1,d}$, і при цьому має місце нерівність
$$
\|f\|_{B^{\boldsymbol{\rho}}_{p',\theta'}}\leqslant C\|f\|_{B^{\boldsymbol{r}}_{p,\theta}},
$$
де $C$~---  деяка константа, яка не залежить від $f$.

\rm \vskip 1 mm

Наведемо ще одне твердження для цілих функцій експоненціального типу, яке одержане С.\,М.~Нікольським~\cite{Nikolsky_1951}, (див., також, \cite[c.~150]{Nikolsky_1969_book}).

\vskip 1 mm
\bf Теорема Б. \it  Якщо
$1\leqslant p_{_1} \leqslant p_{_2} \leqslant \infty$, то для цілої функції
експоненціального типу $g=g_{\boldsymbol{\nu}}\in L_{p_1}(\mathbb{R}^d)$ має місце
``нерівність різних метрик''
\begin{equation}\label{Riz_Metric}
 \|g_{\boldsymbol{\nu}}\|_{L_{p_{_2}}(\mathbb{R}^d)}\leqslant 2^d\left( \prod \limits_{j=1}^d
 \nu_k\right)^{\frac{1}{p_{_1}}-\frac{1}{p_{_2}}}\|g_{\boldsymbol{\nu}}\|_{L_{p_{_1}}(\mathbb{R}^d)}.
\end{equation}\rm \vskip 1 mm

Далі для функції $f\in L_{\infty}(\mathbb{R}^d)$ розглянемо величину
\begin{equation}\label{En_Fourier}
\mathcal{E}_{D_{\boldsymbol{a}^n}}(f)_{\infty}=\|f-S_{\boldsymbol{a}^{n-1}}(f)\|_{\infty}, \ \ n\in \mathbb{N},
\end{equation}
яка називається наближенням функції $f$ її $\boldsymbol{a}^n$-відрізками інтеграла Фур'є.

Відповідно для функціонального класу $F\subset L_{{\infty}}(\mathbb{R}^d)$ покладемо
\begin{equation}\label{En_Fourier_Class}
\mathcal{E}_{D_{\boldsymbol{a}^n}}(F)_{\infty}=\sup\limits_{f\in F} \mathcal{E}_{D_{\boldsymbol{a}^n}}(f)_{\infty}.
\end{equation}

Справедливе таке твердження.

\vskip 1 mm \bf Теорема. \it Нехай $1<p <\infty$, $1\leqslant \theta
\leqslant \infty$. Тоді для  $g(\boldsymbol{r})>\frac{d}{p}$  має місце порядкове співвідношення
\begin{equation} \label{teor_E_inf_aniz}
   {\mathcal{E}}_{D_{\boldsymbol{a}^n}}(B^{\boldsymbol{r}}_{p,\theta})_{\infty}
   \asymp 2^{-n\left(g(\boldsymbol{r})-\frac{d}{p}\right)},
\end{equation}
де $a_j=2^{g(\boldsymbol{r})/r_j}$, $j=\overline{1,d}$. \rm

\vskip 1 mm

Зауважимо, що виконання умови $g(\boldsymbol{r})>\frac{d}{p}$, згідно з теоремою~А, забезпечує належність функцій $f\in B^{\boldsymbol{r}}_{p,\theta}$ до простору $L_{\infty}(\mathbb{R}^d)$.

{\textbf{\textit{Доведення.}}} Спочатку отримаємо
 в (\ref{teor_E_inf_aniz}) оцінку зверху. Оскільки (див. (\ref{vklad})) ${B^{\boldsymbol{r}}_{p,\theta}\subset
 B^{\boldsymbol{r}}_{p,\infty}}=H^{\boldsymbol{r}}_p$, $1\leqslant \theta <\infty$, то шукану оцінку достатньо
 отримати для величини $\mathcal{E}_{D_{\boldsymbol{a}^n}}(H_p^{\boldsymbol{r}})_{\infty}$.

Згідно з  (\ref{f-norm-dek-inf}) для $f\in H^{\boldsymbol{r}}_p$ маємо  ${\|f_{\boldsymbol{a}^s}\|_p\ll 2^{-sg(\boldsymbol{r})}}$. Тому, скориставшись нерівністю Мінковського, нерівністю (\ref{Riz_Metric}), врахувавши, що $a_j=2^{g(\boldsymbol{r})/r_j}$ та беручи до уваги (\ref{g(r)}), отримаємо
$$
{\mathcal{E}}_{D_{\boldsymbol{a}^n}}(f)_{\infty}=\|f-S_{\boldsymbol{a}^{n-1}}(f)\|_{\infty} = \Big\|\sum
 \limits_{s=0}^{\infty}f_{\boldsymbol{a}^s}-S_{\boldsymbol{a}^{n-1}}(f)\Big\|_{\infty}=
$$
 $$
   = \Big\|\sum
 \limits_{s=n}^{\infty}f_{\boldsymbol{a}^s}\Big\|_{\infty}\leqslant \sum
 \limits_{s=n}^{\infty}\|f_{\boldsymbol{a}^s}\|_{\infty} \leqslant \sum
 \limits_{s=n}^{\infty}2^d \left(\prod\limits^d_{j=1} a^s_j \right)^{\frac{1}{p}} \|f_{\boldsymbol{a}^s}\|_{p}\ll
 $$
 $$
 \ll \sum
 \limits_{s=n}^{\infty} 2^{\frac{sd}{p}} \|f_{\boldsymbol{a}^s}\|_{p} \ll
 \sum\limits_{s=n}^{\infty} 2^{\frac{sd}{p}} 2^{- s g(\boldsymbol{r})}= \sum\limits_{s=n}^{\infty}2^{-s\left(g(\boldsymbol{r})-\frac{d}{p}\right)}
 \ll 2^{-n\left(g(\boldsymbol{r})-\frac{d}{p}\right)}.
 $$

%$$
% 2^d \left(\prod\limits^d_{j=1} a^s_j \right)^{\frac{1}{p}}=2^d \left( 2^{\sum\limits_{j=1}^d sg(\boldsymbol{r})/ r_j}\right)^{\frac{1}{p}}\asymp 2^{sd/p}
%&$

Оцінку зверху для величини $\mathcal{E}_{D_{\boldsymbol{a}^n}}(H_p^{\boldsymbol{r}})_{\infty}$ встановлено.

Отримаємо тепер в (\ref{teor_E_inf_aniz}) оцінку знизу. Оскільки ${B^{\boldsymbol{r}}_{p,1}\subset B^{\boldsymbol{r}}_{p,\theta}}$,
$1<\theta\leqslant\infty$, то шукану оцінку  достатньо отримати для величини
$\mathcal{E}_{D_{\boldsymbol{a}^n}}(B^{\boldsymbol{r}}_{p,1})_{\infty}$. Іншими словами достатньо оцінити знизу
величину ${\|f-S_{\boldsymbol{a}^{n-1}}(f)\|_{\infty}}$ для деякої функції $f\in
B^{\boldsymbol{r}}_{p,1}$.

З цією метою розглянемо функцію (див. \cite{Yanchenko_Arxiv+UMG2017})
$$
g_1(\boldsymbol{x})=C_1 2^{-n\left(g(\boldsymbol{r})+\frac{d}{p'}\right)}F_{\boldsymbol{n}}(\boldsymbol{x}),
$$
де $\boldsymbol{n}=(n,\ldots,n)\in \mathbb{N}^d$, $1/p+ 1/p'=1$, $C_1>0$ і
$$
 F_{\boldsymbol{n}}(\boldsymbol{x})=\prod\limits_{j=1}^{d}\sqrt{\frac {2}{\pi}}\ \frac{ \sin{a_j^n x_j}}{x_j}- \prod\limits_{j=1}^{d}\sqrt{\frac {2}{\pi}}\ \frac{ \sin{a_j^{n-1} x_j}}{x_j}
 $$
 та
 $$
F_{0}(\boldsymbol{x})=\prod\limits_{j=1}^{d}\sqrt{\frac {2}{\pi}}\ \frac{ \sin{x_j}}{x_j}.
 $$

 Для перетворення Фур'є функції $F_{\boldsymbol{n}}(\boldsymbol{x})$
 справедливе  співвідношення (див., наприклад, \cite{WangHeping_SunYongsheng_1995})
 $$
 \mathfrak{F}F_{\boldsymbol{n}}(\boldsymbol{x})=\chi_{\boldsymbol{n}}(\boldsymbol{\lambda})=\prod \limits_{j=1}^d \chi_{n}(\lambda_j),
 $$
де

 \begin{minipage}{9 cm}
    $$
      \chi_{n}(\lambda_j)=
 \begin{cases}
    1, & a_j^{n-1}<|\lambda_j|<a_j^{n}, \\
    \frac{1}{2}, & |\lambda_j|=a_j^{n-1} \ \mbox{або} \ |\lambda_j|=a_j^{n}, \\
    0 & \mbox{--- в інших випадках},
 \end{cases}
    $$
\end{minipage}
\begin{minipage}{6 cm}
  $$
   \chi_{0}(x_j)=
 \begin{cases}
    1, & |\lambda_j|<1; \\
    \frac{1}{2}, & |\lambda_j|=1; \\
    0, & |\lambda_j|>1.
 \end{cases}
  $$
\end{minipage}
\vskip 1mm

Відповідно для оберненого перетворення будемо мати
 $$
 \mathfrak{F}^{-1}\chi_{\boldsymbol{n}}(\boldsymbol{\lambda})=F_{\boldsymbol{n}}(\boldsymbol{x}).
 $$

Зазначимо, що $F_{\boldsymbol{n}}(\boldsymbol{x})$~--- ціла функція з $L_p(\mathbb{R}^d)$, носій перетворення Фур'є якої зосереджений в $\Gamma_{\boldsymbol{a}^{\boldsymbol{n}}}$ і крім цього
\begin{equation}\label{F_k_norm_os}
\|F_{\boldsymbol{n}}\|_p\asymp 2^{\frac{dn}{p'}}, \ 1<p<\infty.
\end{equation}

В \cite{Yanchenko_Arxiv+UMG2017} показано, що з деякою константою $C_1>0$ функція $g_1$ належить класу $B^{\boldsymbol{r}}_{p,1}(\mathbb{R}^d)$, а саме,  згідно з (\ref{f-norm-dek}) та (\ref{F_k_norm_os}), маємо
$$
\|g_1\|_{B^{\boldsymbol{r}}_{p,1}}\asymp \sum\limits_{s} 2^{sg(\boldsymbol{r})} \|f_{\boldsymbol{a}^s}(g_1)\|_p\asymp
$$
$$
\asymp  \sum\limits_{s} 2^{sg(\boldsymbol{r})} 2^{-n\left(g(\boldsymbol{r})+\frac{d}{p'}\right)}\|F_{\boldsymbol{n}}\|_p \ll 2^{-n\left(g(\boldsymbol{r})+\frac{d}{p'}\right)} 2^{ng(\boldsymbol{r})} 2^{\frac{dn}{p'}} =1.
$$

Перш ніж перейти до встановлення оцінки знизу в (\ref{teor_E_inf_aniz}),  одержимо порядок величини
\begin{equation}\label{F_k_norm-inf}
\|F_{\boldsymbol{n}}\|_{\infty}=\bigg\|\prod\limits_{j=1}^{d}\sqrt{\frac {2}{\pi}}\ \frac{ \sin{a_j^n x_j}}{x_j}- \prod\limits_{j=1}^{d}\sqrt{\frac {2}{\pi}}\ \frac{ \sin{a_j^{n-1} x_j}}{x_j}\bigg\|_{\infty}.
\end{equation}

Для оцінки зверху будемо мати
$$
\|F_{\boldsymbol{n}}\|_{\infty}= \bigg\|\prod\limits_{j=1}^{d}\sqrt{\frac {2}{\pi}}\ \frac{ \sin{a_j^n x_j}}{x_j} - \prod\limits_{j=1}^{d}\sqrt{\frac {2}{\pi}}\ \frac{ \sin{a_j^{n-1} x_j}}{x_j} \bigg\|_{\infty}\leqslant
$$
$$
\leqslant \bigg\|\prod\limits_{j=1}^{d}\sqrt{\frac {2}{\pi}}\ \frac{ \sin{a_j^n x_j}}{x_j}\bigg\|_{\infty} + \bigg\|\prod\limits_{j=1}^{d}\sqrt{\frac {2}{\pi}}\ \frac{ \sin{a_j^{n-1} x_j}}{x_j} \bigg\|_{\infty} =
$$
$$
 = \sup\limits_{\boldsymbol{x}\in \mathbb{R}^d} \Big|\prod\limits_{j=1}^{d}\sqrt{\frac {2}{\pi}}\ \frac{ \sin{a_j^n x_j}}{x_j}\Big| + \sup\limits_{\boldsymbol{x}\in\mathbb{R}^d} \Big| \prod\limits_{j=1}^{d}\sqrt{\frac {2}{\pi}}\ \frac{ \sin{a_j^{n-1} x_j}}{x_j}\Big| =
$$
$$
 = \left(\frac {2}{\pi}\right)^{\frac{d}{2}}\prod\limits_{j=1}^{d}\sup\limits_{x_j\in\mathbb{R}} \Big| \frac{ \sin{a_j^n x_j}}{x_j}\Big| + \left(\frac {2}{\pi}\right)^{\frac{d}{2}}\prod\limits_{j=1}^{d}\sup\limits_{x_j\in \mathbb{R}} \Big| \frac{ \sin{a_j^{n-1} x_j}}{x_j}\Big|\ll
$$
$$
\ll \left(\prod\limits_{j=1}^{d} a_j^{n}+\prod\limits_{j=1}^{d} a_j^{n-1}\right) = \left(\prod\limits_{j=1}^{d} 2^{ng(\boldsymbol{r})/r_j}+\prod\limits_{j=1}^{d} 2^{(n-1)g(\boldsymbol{r})/r_j}\right)=
$$
\begin{equation}\label{F_k_norm_inf1}
 = \left(2^{dn}+2^{d(n-1)}\right)\ll 2^{dn}.
\end{equation}

Оцінюючи $\|F_{\boldsymbol{n}}\|_{\infty}$ знизу, одержимо
$$
\|F_{\boldsymbol{n}}\|_{\infty} = \bigg\|\prod\limits_{j=1}^{d}\sqrt{\frac {2}{\pi}}\ \frac{ \sin{a_j^n x_j}}{x_j} - \prod\limits_{j=1}^{d}\sqrt{\frac {2}{\pi}}\ \frac{ \sin{a_j^{n-1} x_j}}{x_j} \bigg\|_{\infty}\geqslant
$$
$$
\geqslant \bigg| \Big\|\prod\limits_{j=1}^{d}\sqrt{\frac {2}{\pi}}\ \frac{ \sin{a_j^n x_j}}{x_j}\Big\|_{\infty} - \Big\|\prod\limits_{j=1}^{d}\sqrt{\frac {2}{\pi}}\ \frac{ \sin{a_j^{n-1} x_j}}{x_j} \Big\|_{\infty} \bigg| \gg
$$
\begin{equation}\label{F_k_norm_inf2}
\gg \left|\prod\limits_{j=1}^{d} a_j^{n} - \prod\limits_{j=1}^{d} a_j^{n-1}\right| \gg (2^{dn}- 2^{d(n-1)})\gg 2^{dn}.
\end{equation}

Співставивши (\ref{F_k_norm_inf1}) і (\ref{F_k_norm_inf2}), можемо записати порядкове співвідношення
\begin{equation}\label{F_k_norm-inf-os}
\|F_{\boldsymbol{n}}\|_{\infty}\asymp 2^{dn}.
\end{equation}

Оскільки, для функції $g_1$  має місце співвідношення $S_{\boldsymbol{a}^{n-1}}(g_1)=0$, то скориставшись (\ref{F_k_norm-inf-os}) приходимо до шуканої оцінки знизу
$$
\mathcal{E}_{D_{\boldsymbol{a}^{n}}}(B^{\boldsymbol{r}}_{p,1})_{\infty}\geqslant
\mathcal{E}_{D_{\boldsymbol{a}^{n}}}(g_1)_{\infty}=\|g_1-S_{\boldsymbol{a}^{n-1}}(g_1)\|_{\infty}= \|g_1\|_{\infty}\gg
$$
$$
\gg2^{-n\left(g(\boldsymbol{r})+\frac{d}{p'}\right)}\|F_{\boldsymbol{n}}\|_{\infty}\gg 2^{-n\left(g(\boldsymbol{r})+\frac{d}{p'}\right)} 2^{dn}=2^{-n\left(g(\boldsymbol{r})-\frac{d}{p}\right)}.
$$

Оцінку знизу в  (\ref{teor_E_inf_aniz}) встановлено. Теорему доведено.

\vskip 1 mm

  На завершення роботи зробимо деякі коментарі.

  Точні за порядком значення величини $\mathcal{E}_{D_{\boldsymbol{a}^{n}}}(B^{\boldsymbol{r}}_{p,\theta})_{q}$, $1<p\leqslant q < \infty$ встановлено в \cite{Yanchenko_Arxiv+UMG2017}.

  У випадку $r_1=\ldots=r_d=r$, тобто для ізотропних класів Нікольського--Бєсова $B^r_{p,\theta}(\mathbb{R}^d)$, оцінку (\ref{teor_E_inf_aniz}) встановлено у роботі~\cite{Yanchenko_2015UMG}. В одновимірному  випадку $(d=1)$ анізотропні класи  Нікольського--Бєсова збігаються з класами Нікольського--Бєсова мішаної гладкості, які досліджувалися у роботах  \cite{WangHeping_SunYongsheng_1995}, \cite{Yanchenko_YMG_2010_8}.

Зазначимо також, що деякі  апроксимативні характеристик
ізотропних та анізотропних класів Нікольського--Бєсова періодичних функцій багатьох змінних досліджувалися, зокрема, у роботах~\cite{Romanyuk_UMG_2009}\,--\,\cite{Myronyuk_UMG_2016_8}.

\vskip 3 mm

\textbf{Contact information:}
Department of the Theory of Functions, Institute of Mathematics of National
Academy of Sciences of Ukraine, 3, Tereshenkivska st., 01004, Kyiv, Ukraine.

\vskip 3 mm

E-mail: \href{mailto:Yan.Sergiy@gmail.com}{Yan.Sergiy@gmail.com}

\end{document}